\documentclass[oneside,english]{amsart}
\usepackage[T1]{fontenc}
\usepackage[latin1]{inputenc}
\usepackage{amssymb}

\makeatletter

\providecommand{\LyX}{L\kern-.1667em\lower.25em\hbox{Y}\kern-.125emX\@}

 \theoremstyle{plain}    
 \newtheorem{thm}{Theorem}[section]
 \numberwithin{equation}{section} 
 \numberwithin{figure}{section} 
 \newcommand{\lyxaddress}[1]{
   \par {\raggedright #1 
   \vspace{1.4em}
   \noindent\par}
 }
 \theoremstyle{plain}    
 \newtheorem{cor}[thm]{Corollary} 
 \theoremstyle{plain}    
 \newtheorem{lem}[thm]{Lemma} 
 \theoremstyle{plain}    
 \newtheorem{prop}[thm]{Proposition} 

\usepackage{babel}
\makeatother
\begin{document}

\title{Multiplicative Monotonic Convolution}

\author{Hari Bercovici}

\thanks{The author was supported in part by a grant from the National Science
Foundation.}

\begin{abstract}
We show that the monotonic independence introduced by Muraki can also
be used to define a multiplicative convolution. We also find a method
for the calculation of this convolution based on an appropriate form
of the Cauchy transform. We also discuss infinite divisibility in
the multiplicative monotonic context.
\end{abstract}
\maketitle

\section{Introduction}

Consider an algebraic probability space, that is, a pair $(\mathfrak{A},\varphi )$
where $\mathfrak{A}$ is a unital complex algebra, and $\varphi :\mathfrak{A}\to \Bbb C$
is a linear functional satisfying $\varphi (1)=1$. Muraki {[}5{]}
introduced the concept of \emph{monotonic independence} for elements
of $\mathfrak{A}$, which we will now review. Let $\mathfrak{A}_{1},\mathfrak{A}_{2}$
be two subalgebras of $\mathfrak{A}$; it is not assumed that either
of these subalgebras contains the unit. These algebras are said to
be \emph{monotonically independent} if the following two conditions
are satisfied:

\begin{enumerate}
\item for every $x_{1},y_{1}\in \mathfrak{A}_{1}$ and $x_{2}\in \mathfrak{A}_{2}$,
we have $x_{1}x_{2}y_{1}=\varphi (x_{2})x_{1}y_{1}$;
\item for every  $x_{1}\in \mathfrak{A}_{1}$ and $x_{2},y_{2}\in \mathfrak{A}_{2}$,
we have $\varphi (x_{2}x_{1}y_{2})=\varphi (x_{2})\varphi (x_{1})\varphi (y_{2})$,
$\varphi (x_{2}x_{1})=\varphi (x_{2})\varphi (x_{1})$, and $\varphi (x_{1}y_{2})=\varphi (x_{1})\varphi (y_{2})$.
\end{enumerate}
Proceeding inductively, the algebras $\mathfrak{A}_{1},\mathfrak{A}_{2},\dots ,\mathfrak{A}_{n}$
are said to be monotonically independent if $\mathfrak{A}_{1},\mathfrak{A}_{2},\dots ,\mathfrak{A}_{n-1}$
are monotonically independent, and the algebras $\mathfrak{A}',\mathfrak{A}_{n}$
are monotonically independent, where $\mathfrak{A}'$ is the (generally
nonunital) algebra generated by $\mathfrak{A}_{1}\cup \mathfrak{A}_{2}\cup \cdots \cup \mathfrak{A}_{n-1}$.
More generally, if $I$ is a totally ordered set, and $(\mathfrak{A}_{i})_{i\in I}$
is a family of subalgebras of $\mathfrak{A}$, this family is said
to be monotonically independent if the algebras $\mathfrak{A}_{i_{1}},\mathfrak{A}_{i_{2}},\dots ,\mathfrak{A}_{i_{n}}$
are monotonically independent for any choice of indices $i_{1}<i_{2}<\cdots <i_{n}$.
A family $(x_{i})_{i\in I}$ of elements of $\mathfrak{A}$ is said
to be monotonically independent if the (generally nonunital) subalgebras
$\mathfrak{A}_{i}$ generated by $x_{i}$ form a monotonically independent
family.

The distribution $\mu _{x}$ of an element $x$ of $\mathfrak{A}$
(a.k.a. a random variable) is the linear functional defined on the
polynomial algebra $\Bbb C[X]$ by the formula \[
\mu _{x}(p)=\varphi (p(x)),\quad p\in \Bbb C[X].\]
 Clearly, $\mu _{x}$ is entirely determined by the sequence $\mu _{x}(X^{n})=\varphi (x^{n})$
of moments of $x.$ A functional $\mu $ on $\Bbb C[X]$ is the distribuion
of some random variable if and only if $\mu (1)=1$. The set of these
functionals, endowed with the weak{*} topology, will be denoted $\mathfrak{M}$.

Muraki {[}6{]} observed that, given monotonically independent random
variables $x_{1,}x_{2}$, the distribution of $x_{1}+x_{2}$ only
depends on $\mu _{x_{1}},\mu _{x_{2}}$. This gives rise to a binary
operation $\triangleright $ on $\mathfrak{M}$, called monotonic
convolution. It was also shown in {[}6{]} how to calculate monotonic
convolutions using moment generating functions. 

It is also true that the distribution of $x_{1}x_{2}$ only depends
on $\mu _{x_{1}},\mu _{x_{2}}$ if $x_{1},x_{2}$ are monotonically
independent, but the dependence is rather trivial. Indeed, if $n\ge 1$,
property (1) above yields\[
(x_{1}x_{2})^{n}=\varphi (x_{2})^{n-1}x_{1}^{n}x_{2},\]
and then from property (2)\[
\varphi ((x_{1}x_{2})^{n})=\varphi (x_{2})^{n}\varphi (x_{1}^{n}).\]
In other words, the product $x_{1}x_{2}$ has the same distribution
as $\alpha x_{1}$, with $\alpha =\varphi (x_{2})$.

A more interesting result is obtained by considering variables $x_{1,}x_{2}\in \mathfrak{A}$
such that the centered variables $x_{1}-c_{1},x_{2}-c_{2}$ are monotonically
independent, where $c_{1},c_{2}$ are scalars. It is again easy to
see that, under this condition, $\mu _{x_{1}x_{2}}$ depends only
on the distributions of $x_{1},x_{2}$ and on the numbers $c_{1},c_{2}$.
This yields a new operation $\circlearrowright $ on $\mathfrak{M}\times \Bbb C$,
called \emph{multiplicative monotonic convolution}, such that\[
(\mu _{x_{1}x_{2}},c_{1}c_{2})=(\mu _{x_{1}},c_{1})\circlearrowright (\mu _{x_{2}},c_{2})\]
 if $x_{1}-c_{1}$ and $x_{2}-c_{2}$ are monotonically independent.
It is interesting to note that, under this condition, $\mu _{x_{1}x_{2}}=\mu _{x_{2}x_{1}}$,
but the operation $\circlearrowright $ is not commutative, since
monotonic independence itself is not a symmetric relation. While $\circlearrowright $
is not an operation on $\mathfrak{M}$ itself, there are two ways
in which it induces such an operation. The first one is obtained identifying
$\mathfrak{M}$ with the subset $\{(\mu ,1):\mu \in \mathfrak{M}\}$;
we will use the same notation for the operation induced this way,
that is\[
(\mu _{1},1)\circlearrowright (\mu _{2},1)=(\mu _{1}\circlearrowright \mu _{2},1),\quad \mu _{1},\mu _{2}\in \mathfrak{M}.\]
The second one is obtained identifying $\mathfrak{M}$ with the subset
$\{(\mu ,\mu (X)):\mu \in \mathfrak{M}$; we will use the notation
$\circlearrowright _{0}$ for this operation, so that\[
(\mu _{1},\mu _{1}(X))\circlearrowright (\mu _{2},\mu _{2}(X))=(\mu _{1}\circlearrowright _{0}\mu _{2},\mu _{1}(X)\mu _{2}(X)),\quad \mu _{1},\mu _{2}\in \mathfrak{M}.\]
The operation $\circlearrowright $ on $\mathfrak{M}$ has the advantage
that it is easily extended to measures with unbounded supports. On
the other hand, $\circlearrowright _{0}$ has the advantage that convolution
with a Dirac point mass has the natural dilation effect.

We will show that multiplicative monotonic convolution can be calculated
in terms of an appropriate moment generating series. We will deduce
from this that the multiplicative monotonic convolution of two probability
measures on the unit circle is again a probability measure on the
unit circle. Analogously, the multiplicative monotonic convolution
of two compactly supported probability measures on $\Bbb R_{+}=[0,+\infty )$
is a measure of the same kind. As mentioned above, the operation $\circlearrowright $
extends to arbitrary probability measures on $\Bbb R_{+}$. It is
not clear whether the same is true for $\circlearrowright _{0}$.multiplicative
monotonic convolution can be extended to arbitrary probability measures
on $\Bbb R_{+}$. In the case of probability measures on $\Bbb R_{+}$
and $\Bbb T$ we will give a description of one-parameter convolution
semigroups and of infinitely divisible measures, at least for compact
supports. This was done by Muraki {[}6{]} for additive monotonic convolution
semigroups of compactly supported measures on $\Bbb R$.

Our approach in calculating multiplicative monotonic convolutions
is related to the one we used in {[}2{]} to approach additive monotonic
convolution, rather than the original combinatorial approach of {[}6{]}.

\section{Realization of Monotonically Independent Variables}

For every distribution $\mu \in \mathfrak{M}$ we will consider the
formal power series\[
\psi _{\mu }(z)=\sum _{n=1}^{\infty }\mu (X^{n})z^{n},\quad \eta _{\mu }(z)=\frac{\psi _{\mu }(z)}{1+\psi _{\mu }(z)}.\]
 If $x$ is a random variable, we will also use the notation $\psi _{x}=\psi _{\mu _{x}},\eta _{x}=\eta _{\mu _{x}}$.
The calculation of multiplicative monotonic convolution will involve
the series $\eta _{\mu }$. 

Consider a Hilbert space $\mathfrak{H}$, and a unit vector $\xi \in \mathfrak{H}$.
The algebra $\mathfrak{L}(\mathfrak{H})$ of bounded linear operators
on $\mathfrak{H}$ becomes a probability space with the vector functional
$\varphi _{\xi }(x)=(x\xi ,\xi )$, $x\in \mathfrak{L}(\mathfrak{H})$.
Consider now the space $\mathfrak{H}'=\mathfrak{H}\otimes \mathfrak{H}$
and the unit vector $\xi '=\xi \otimes \xi $. Following Muraki {[}6{]},
one produces two monotonically independent copies $\mathfrak{A}_{1},\mathfrak{A}_{2}$
of $\mathfrak{L}(\mathfrak{H})$ in $(\mathfrak{L}(\mathfrak{H}'),\varphi _{\xi '})$
as follows: \[
\mathfrak{A}_{1}=\{x\otimes p:x\in \mathfrak{L}(\mathfrak{H})\},\quad \mathfrak{A}_{2}=\{1\otimes x:x\in \mathfrak{L}(\mathfrak{H})\},\]
 where $p$ denotes the rank one projection onto the space generated
by $\xi $, and $1$ denotes the identity operator on $\mathfrak{H}$.

For $x\in \mathfrak{L}(\mathfrak{H})$, the formal power series $\psi _{x}(z),\eta _{x}(z)$
are actually convergent, at least for $|z|<1/\| x^{-1}\| $. Assume
now that $\mathfrak{H}$ has an orthonormal basis $(\xi _{j})_{j=1}^{\infty }$,
and $\xi =\xi _{0}$. Consider the shift $s\in \mathfrak{L}(\mathfrak{H})$
defined by $s\xi _{j}=\xi _{j+1}$ for all $j$. We will be interested
in elements $x\in \mathfrak{L}(\mathfrak{H})$ of the form $x=(1+s)u(s^{*})$,
where $u\in \Bbb C[X]$ is a polynomial. It is easy to see that the
distributions $\mu _{x}$ of these operators form a dense subset in
$\mathfrak{M}$. Moreover, as shown by Haagerup (see Theorem 2.3.(a)
in {[}4{]}), the generating function $\psi _{x}$ is easily related
to $u$.

\begin{lem}
If $x=(1+s)u(s^{*})$, where $u$ is a polynomial with $u(0)\ne 0$,
then\[
\psi _{x}\left(\frac{z}{(1+z)u(z)}\right)=z\]
for sufficintly small $|z|$.
\end{lem}
We can now state the main result of this section.

\begin{thm}
Consider two distributions $\mu _{1},\mu _{2}\in \mathfrak{M}$, constants
$c_{1},c_{2}\in \Bbb C$ and the multiplicative monotonic convolution
$(\mu ,c_{1}c_{2})=(\mu _{1},c_{1})\circlearrowright (\mu _{2},c_{2})$.
We have then\[
\eta _{\mu }(z)=\eta _{\mu _{1}}\left(\frac{1}{c_{1}}\eta _{\mu _{2}}(c_{1}z)\right),\]
if $c_{1}\ne 0$, and\[
\eta _{\mu }(z)=\eta _{\mu _{1}}(\eta '_{\mu _{2}}(0)z),\]
if $c_{1}=0.$
\end{thm}
\begin{proof}
Clearly the operation $\circlearrowright $ is continuous, and therefore
it will suffice to prove the theorem for $\mu _{1},\mu _{2}$ in a
dense family of distributions, for instance the family of distributions
obtained from the random variables $(1+s)u(s^{*})$, where $u$ is
a polynomial with $u(0)\ne 0$. Assume then that $u_{1},u_{2}$ are
two polynomials which do not vanish at the origin. We consider the
variables $x_{1},x_{2}$ in ($\mathfrak{L}(\mathfrak{H}\otimes \mathfrak{H}),\varphi _{\xi _{0}\otimes \xi _{0}})$
defined by\[
x_{1}=c_{1}\otimes (1-p)+(1+s)u_{1}(s^{*})\otimes p,\quad x_{2}=(1+1\otimes s)u_{2}(1\otimes s^{*})=1\otimes [(1+s)u_{2}(s^{*})].\]
These variables have the required property that $x_{1}-c_{1},x_{2}-c_{2}$
belong to monotonically independent subalgebras. Moreover, it is easy
to see that $x_{1},x_{2}$ have the same distributions as $(1+s)u_{1}(s^{*}),(1+s)u_{2}(s^{*})$,
so that\[
\psi _{x_{1}}\left(\frac{z}{(1+z)u_{1}(z)}\right)=z,\quad \psi _{x_{2}}\left(\frac{z}{(1+z)u_{2}(z)}\right)=z\]
for sufficiently small $z$.

Consider now the vectors $\xi _{\lambda }=\sum _{n=0}^{\infty }\lambda ^{n}\xi _{n}\in \mathfrak{H}$
defined for $|\lambda |<1$. Since\[
s\xi _{\mu }=\frac{1}{\mu }(\xi _{\mu }-\xi _{0}),\quad u(s^{*})\xi _{\mu }=u(\mu )\xi _{\mu }\]
for $\mu \ne 0$, we can easily calculate\[
x_{2}(\xi _{\lambda }\otimes \xi _{\mu })=(1\otimes (1+s))u_{2}(\mu )\xi _{\lambda }\otimes \xi _{\mu }=u_{2}(\mu )\xi _{\lambda }\otimes \left(\xi _{\mu }+\frac{1}{\mu }(\xi _{\mu }-\xi _{0})\right).\]
 Then we obtain for $\lambda \ne 0\ne \mu $\begin{eqnarray*}
 & x_{1}x_{2}(\xi _{\lambda }\otimes \xi _{\mu })=c_{1}u_{2}(\mu )\xi _{\lambda }\otimes \left(\left (1+\frac{1}{\mu }\right )(\xi _{\mu }-\xi _{0})\right)+u_{2}(\mu )((1+s)\otimes p)\left[u_{1}(\lambda )\xi _{\lambda }\otimes \xi _{0}\right] & \\
 & =u_{2}(\mu )\left[c_{1}\left(1+\frac{1}{\mu }\right)\xi _{\lambda }\otimes (\xi _{\mu }-\xi _{0})+u_{1}(\lambda )\xi _{\lambda }\otimes \xi _{0}+\frac{u_{1}(\lambda )}{\lambda }(\xi _{\lambda }-\xi _{0})\otimes \xi _{0}\right] & \\
 & =u_{2}(\mu )\left[c_{1}\left(1+\frac{1}{\mu }\right)\xi _{\lambda }\otimes \xi _{\mu }+\left(u_{1}(\lambda )\left(1+\frac{1}{\lambda }\right)-c_{1}\left(1+\frac{1}{\mu }\right)\right)\xi _{\lambda }\otimes \xi _{0}-\frac{u_{1}(\lambda )}{\lambda }\xi _{0}\otimes \xi _{0}\right]. & 
\end{eqnarray*}
This equation can be simplified when\[
u_{1}(\lambda )\left(1+\frac{1}{\lambda }\right)-c_{1}\left(1+\frac{1}{\mu }\right)=0,\]
in which case it becomes\[
x_{1}x_{2}(\xi _{\lambda }\otimes \xi _{\mu })=\frac{1}{z}(\xi _{\lambda }\otimes \xi _{\mu })-\frac{u_{1}(\lambda )u_{2}(\mu )}{\lambda }\xi _{0}\otimes \xi _{0},\quad \text {with}\quad \frac{1}{z}=c_{1}u_{2}(\mu )\left(1+\frac{1}{\mu }\right).\]
This can then be rewritten as\[
(1-zx_{1}x_{2})^{-1}\xi _{0}\otimes \xi _{0}=\frac{\lambda }{zu_{1}(\lambda )u_{2}(\mu )}\xi _{\lambda }\otimes \xi _{\mu },\]
so that \[
\varphi ((1-zx_{1}x_{2})^{-1})=((1-zx_{1}x_{2})^{-1}\xi _{0}\otimes \xi _{0},\xi _{0}\otimes \xi _{0})=\frac{\lambda }{zu_{1}(\lambda )u_{2}(\mu )}.\]
The constant in the right-hand side of this equation is now easily
calculated:\[
\frac{\lambda }{zu_{1}(\lambda )u_{2}(\mu )}=\frac{\lambda c_{1}(1+1/\mu )}{u_{1}(\lambda )}=\lambda \left(1+\frac{1}{\lambda }\right)=\lambda +1,\]
yielding then\[
\psi _{x_{1}x_{2}}(z)=\varphi _{\xi _{0}\otimes \xi _{0}}((1-zx_{1}x_{2})^{-1})-1=\lambda .\]
These calculations hold for $|\lambda |\ne 0$ sufficiently small,
because the associated numbers $\mu $ and $z$ are also small, and
$\mu \ne 0$. Observe now that the identity \[
\frac{1}{c_{1}z}=u_{2}(\mu )\left(1+\frac{1}{\mu }\right)\]
means that $\psi _{x_{2}}(c_{1}z)=\mu ,$ while \[
u_{1}(\lambda )\left(1+\frac{1}{\lambda }\right)=c_{1}\left(1+\frac{1}{\mu }\right)\]
means that\[
\lambda =\psi _{x_{1}}\left(\frac{1}{c_{1}(1+1/\mu )}\right).\]
Combining these identities we see that\[
\psi _{x_{1}x_{2}}(z)=\psi _{x_{1}}\left(\frac{1}{c_{1}(1+1/\mu )}\right)=\psi _{x_{1}}\left(\frac{1}{c_{1}}\frac{\psi _{x_{2}}(c_{1}z)}{1+\psi _{x_{2}}(c_{1}z)}\right)=\psi _{x_{1}}\left(\frac{1}{c_{1}}\eta _{x_{2}}(c_{1}z)\right).\]
The identity above shows that \[
\eta _{x_{1}x_{2}}(z)=\eta _{x_{1}}\left(\frac{1}{c_{1}}\eta _{x_{2}}(c_{1}z)\right)\]
for uncountably many values of $z$. We deduce that the identity in
the statement holds in the generic particular case $\mu _{1}=\mu _{x_{1}},\mu _{2}=\mu _{x_{2}}$.
\end{proof}
The two convolutions on $\mathfrak{M}$ are now easily described.

\begin{cor}
Given measures $\mu _{1},\mu _{2}\in \mathfrak{M}$, we have\[
\eta _{\mu _{1}\circlearrowright \mu _{2}}(z)=\eta _{\mu _{1}}(\eta _{\mu _{2}}(z)),\]
and\[
\eta _{\mu _{1}\circlearrowright _{0}\mu _{2}}(z)=\eta _{\mu _{1}}\left (\frac{1}{\alpha }\eta _{\mu _{2}}(\alpha z)\right ),\]
with $\alpha =\mu _{1}(X)=\eta _{\mu _{1}}'(0)$. The fraction $\eta _{\mu _{2}}(\alpha z)/\alpha $
must be interpreted as $\eta _{\mu _{2}}'(0)z$ in case $\alpha =0$.
\end{cor}

\section{Measures on the Positive Half-Line}

If $\mu $ is a probability measure on $\Bbb R_{+}$ one can define\[
\psi _{\mu }(z)=\int _{0}^{\infty }\frac{zt}{1-zt}\, d\mu (t),\quad \eta _{\mu }(z)=\frac{\psi _{\mu }(z)}{1+\psi _{\mu }(z)}\]
for every $z\in \Omega =\Bbb C\setminus \Bbb R_{+}$. These functions
are analytic, and moreover $\eta _{\mu }(\Omega )\subset \Omega $.
More precisely,\[
\eta _{\mu }(0-)=0,\quad \eta _{\mu }(\overline{z})=\overline{\eta _{\mu }(z),}\quad \text {and}\quad \pi \ge \arg \eta _{\mu }(z)\ge \arg z,\quad \text {for}\, \, z\in \Omega ,\Im z>0,\]
where $\eta _{\mu }(0-)=\lim _{t\uparrow 0}\eta _{\mu }(t)$. Moreover,
as seen in {[}1{]}, these conditions characterize the functions $\eta _{\mu }$
among all analytic functions defined on $\Omega $. The measure $\mu $
is compactly supported if and only if the function $\eta _{\mu }$
is analytic in a neighborhood of the origin. In this case, $\mu $
is entirely determined by the Taylor coefficients of $\eta _{\mu }$,
and the power series of $\eta _{\mu }$ at zero is precisely the formal
power series denoted by the same symbol in the preceding section,
provided that we view $\mu $ as an element of $\mathfrak{M}$ by
setting $\mu (X^{n})=\int _{0}^{\infty }t^{n}\, d\mu (t)$. We can
thus identify the collection of compactly supported measures on $\Bbb R_{+}$
with a subset of $\mathfrak{M}$.

\begin{prop}
If $\mu _{1},\mu _{2}$ are compactly supported probability measures
on $\Bbb R_{+}$ then both $\mu _{1}\circlearrowright \mu _{2}$ and
$\mu _{1}\circlearrowright _{0}\mu _{2}$ are compactly supported
probability measure on $\Bbb R_{+}$.
\end{prop}
\begin{proof}
If $\mu _{1}=\delta _{0}$ is Dirac measure at zero, then clearly
$\mu _{1}\circlearrowright _{0}\mu _{2}=\delta _{0}$. Otherwise,
the number $\alpha =\eta '_{\mu _{1}}(0)=\int _{0}^{\infty }t\: d\mu (t)$
is different from zero, and therefore\[
\eta _{\mu _{1}\circlearrowright _{0}\mu _{2}}(z)=\eta _{\mu _{1}}\left(\frac{1}{\alpha }\eta _{\mu _{2}}(\alpha z)\right).\]
This shows that $\eta _{\mu _{1}\circlearrowright _{0}\mu _{2}}(z)$
makes sense for every $z\in \Omega $, and it is an analytic function
of the form $\eta _{\mu }$ for some compactly supported probability
measure $\mu $ on $\Bbb R_{+}$. Clearly then $\mu =\mu _{1}\circlearrowright _{0}\mu _{2}$.
The case of $\mu _{1}\circlearrowright \mu _{2}$ is treated similarly.
\end{proof}
There is a different argument for the preceding result, based on the
multiplication of positive random variables. Observe first that the
existence of monotonically independent variables is, generally, incompatible
with the functional linear $\varphi $ being a trace. Indeed, if $\mathfrak{A}_{1},\mathfrak{A}_{2}$
are monotonically independent in $(\mathfrak{A},\varphi )$, and $x_{1}\in \mathfrak{A}_{1},x_{2},y_{2}\in \mathfrak{A}_{2},$
then\[
\varphi (x_{2}x_{1}y_{2})-\varphi (x_{1}y_{2}x_{2})=\varphi (x_{1})[\varphi (x_{2})\varphi (y_{2})-\varphi (x_{2}y_{2})].\]
Thus, if $\varphi $ is a trace, either $\varphi |\mathfrak{A}_{1}$
is identically zero, or $\varphi |\mathfrak{A}_{2}$ is multiplicative.
There is however a remnant of the trace property, for instance when
$\mathfrak{A}_{1}$ is commutative.

\begin{lem}
Assume that $\mathfrak{A}_{1},\mathfrak{A}_{2}$ are monotonically
independent in $(\mathfrak{A},\varphi )$, and $\varphi |\mathfrak{A}_{1}$is
a trace. Then we have $\varphi (xy)=\varphi (yx)$ for any $x$ in
the unital algebra generated by $\mathfrak{A}_{1}$, and any $y$
in the unital algebra generated by $\mathfrak{A}_{1}\cup \mathfrak{A}_{2}$.
\end{lem}
\begin{proof}
Since both sides of the identity to be proved are bilinear in $(x,y)$,
it suffices to prove it when $x\in \mathfrak{A}_{1}$, and $y$ is
a product of elements in $\mathfrak{A}_{1}\cup \mathfrak{A}_{2}$,
with at least one factor in $\mathfrak{A}_{2}$. Thus $y$ has the
form\[
y=x_{1}y_{1}\cdots x_{n}y_{n}x_{n+1},\]
where $n\ge 1,$ $y_{1},y_{2},\dots ,y_{n}\in \mathfrak{A}_{2}$,
$x_{2},\dots ,x_{n}\in \mathfrak{A}_{1}$, and $x_{1},x_{n+1}\in \mathfrak{A}_{1}\cup \{1\}$.
Monotonic independence allows us to calculate\[
\varphi (xy)-\varphi (yx)=[\varphi (xx_{1}x_{2}\cdots x_{n+1})-\varphi (x_{1}x_{2}\cdots x_{n+1}x)]\prod _{j=1}^{n}\varphi (y_{j}),\]
and the conclusion follows because $\varphi |\mathfrak{A}_{1}$ is
a trace.
\end{proof}
\begin{cor}
Let $x_{1},x_{2}$ be two random variables, and $c_{1},c_{2}\in \Bbb C$
be such that $x_{1}-c_{1}$ and $x_{2}-c_{2}$ are monotonically independent.
Then the variables $x_{1}^{2}x_{2},x_{1}x_{2}x_{1},$ and $x_{2}x_{1}^{2}$
have the same distribution.
\end{cor}
In particular, if the probability space is $(\mathfrak{L}(\mathfrak{H}),\varphi _{\xi })$,
and $x_{1},x_{2}$ are selfadjoint, it follows that $x_{1}^{2}x_{2}$
has the same distribution as the selfadjoint variable $x_{1}x_{2}x_{1}$.
If $\mu _{1},\mu _{2}$ are compactly supported measures on $\Bbb R_{+}$,
then we can always find random variables $y_{1},y_{2}\in \mathfrak{L}(\mathfrak{H})$
such that $y_{1},y_{2}$ are positive operators, and $\mu _{y_{1}^{2}}=\mu _{1},\mu _{y_{2}}=\mu _{2}$.
We can then define new variables\[
x_{1}=y_{1}\otimes p+c_{1}^{1/2}\otimes (1-p),\quad x_{2}=1\otimes y_{2}\]
in $(\mathfrak{L}(\mathfrak{H}\otimes \mathfrak{H}),\varphi _{\xi \otimes \xi })$
which have the same distributions as $y_{1},y_{2}$, and $x_{1}-c_{1}^{1/2},x_{2}-c_{2}$
are monotonically independent. Considering now $c_{1}=c_{2}=1$ or
$c_{j}=\mu _{j}(X)$, we see that $\mu _{1}\circlearrowright \mu _{2}$,
and respectively $\mu _{1}\circlearrowright _{0}\mu _{2}$, is the
distribution of the positive random variable $x_{1}x_{2}x_{1}$. Moreover,
the inequality $\| x_{1}x_{2}x_{1}\| \le \| x_{1}^{2}\| \| x_{2}\| $,
and the fact that $y_{1},y_{2}$ can be chosen so that the spectra
of $y_{1}^{2},y_{2}$ coincide with the supports of $\mu _{1},\mu _{2}$,
yield the following result.

\begin{cor}
Let $\mu _{1},\mu _{2}$ be probability measures on $\Bbb R_{+}$such
that the support of $\mu _{j}$ is contained in the interval $[\alpha _{j},\beta _{j}]\subset \Bbb R_{+}$,
where $\alpha _{j}\le 1\le \beta _{j}$ for $j=1,2$. Then the supports
of $\mu _{1}\circlearrowright \mu _{2}$ and $\mu _{1}\circlearrowright _{0}\mu _{2}$
are contained in $[\alpha _{1}\alpha _{2},\beta _{1}\beta _{2}]$.
\end{cor}
We will need an inclusion in the opposite direction. In the following
proof we will use the fact that the measure $\mu $ can be recovered
from the imaginary parts of the limits of the function $\psi _{\mu }$
or $\eta _{\mu }$ at points on the real line. The relevant fact is
as follows: if $(a,b)\subset \Bbb R_{+}$ is an open interval such
that $\lim _{\theta \downarrow 0}\arg \eta _{\mu }(re^{i\theta })=0$
for every $r\in (a,b)$ , then $\mu ((1/b,1/a))=0$.

We will denote by $\text {supp}(\mu )$ the supremum of the support
of a measure $\mu $ on $\Bbb R_{+}.$

\begin{prop}
For any compactly supported probability measures $\mu _{1},\mu _{2}$
on $\Bbb R_{+}$, we have \[
\text {supp}(\mu _{2})\subset \text {supp}(\mu _{1}\circlearrowright \mu _{2}),\]
and \[
\left(\int _{0}^{\infty }t\, d\mu _{1}(t)\right)\text {supp}(\mu _{2})\subset \text {supp}(\mu _{1}\circlearrowright _{0}\mu _{2}).\]

\end{prop}
\begin{proof}
We provide the argument for $\mu =\mu _{1}\circlearrowright _{0}\mu _{2}$.
Assume that an interval $(a,b)$ is disjoint from the support of $\mu $,
so that the function $\eta _{\mu }$ is analytic and real-valued on
the interval $(1/b,1/a)$. Now\[
\eta _{\mu }(z)=\eta _{\mu _{1}}\left(\frac{1}{\alpha }\eta _{\mu _{2}}(\alpha z)\right)\]
for $z\notin \Bbb R_{+}$, with\[
\alpha =\eta '_{\mu _{1}}(0)=\int _{0}^{\infty }t\, d\mu _{1}(t).\]
We deduce that\[
\lim _{\theta \downarrow 0}\arg \eta _{\mu _{2}}(\alpha re^{i\theta })=\lim _{\theta \downarrow 0}\arg \frac{1}{\alpha }\eta _{\mu _{2}}(\alpha re^{i\theta })\le \lim _{\theta \downarrow 0}\arg \eta _{\mu }(re^{i\theta })=0\]
for every $r\in (1/b,1/a)$. As noted before the statement of the
proposition, this implies that the support of the measure $\mu _{2}$
contains no points in $(a/\alpha ,b/\alpha )$. In other words,\[
\text {supp}(\mu _{2})\subset \frac{\text {supp}(\mu _{1}\circlearrowright _{0}\mu _{2})}{\int _{0}^{\infty }t\, d\mu _{1}(t)},\]
as claimed.
\end{proof}
It is now fairly easy to find the multiplicative monotonic convolution
semigroups. These are simply families $\{\mu _{\tau }:\tau \ge 0\}$
of compactly supported probability measures on $\Bbb R_{+}$ such
that $\mu _{0}=\delta _{1}$, $\mu _{\tau +\tau '}=\mu _{\tau }\circlearrowright \mu _{\tau '}$
(or $\mu _{\tau +\tau '}=\mu _{\tau }\circlearrowright _{0}\mu _{\tau '}$)
for $\tau ,\tau '\ge 0$, and the map $\tau \mapsto \mu _{\tau }$
is continuous. The topology on probability measures will be the one
inherited from $\mathfrak{M}$, but in this case it is precisely the
topology of weak convergence of probability measures. Indeed, in the
case of $\circlearrowright $-semigroups the support of $\mu _{\tau }$
is contained in the support of $\mu _{1}$ for $\tau \le 1$, and
it is immediate that the map $\tau \mapsto \mu _{\tau }$ is continuous
when we consider the weak topology on the collection of probability
measures. Similarly, in the case of $\circlearrowright _{0}$-semigroups,
observe first that the function $\alpha (\tau )=\int _{0}^{\infty }t\, d\mu _{\tau }(t)$
is continuous, and $\alpha (\tau +\tau ')=\alpha (\tau )\alpha (\tau ')$.
We conclude that $\alpha (\tau )=e^{a\tau }$, with $a=\log \alpha (1)\in \Bbb R$.
The preceding result now shows that the support of $\mu _{\tau }$
is uniformly bounded when $\tau $ runs in a bounded set. Indeed,
we see that\[
\text {supp}(\mu _{\tau })\subset \frac{\text {supp}(\mu _{T})}{e^{a(T-\tau )}}\]
for $\tau \in [0,T]$. It is again easy to conclude that the map $\tau \mapsto \mu _{\tau }$
is continuous when we consider the weak topology on the collection
of probability measures.

\begin{thm}
Consider a $\circlearrowright _{0}$-semigroup $\{\mu _{\tau }:\tau \ge 0\}$
of compactly supported probability measures on $\Bbb R_{+}$, and
let $a\in \Bbb R$ be such that\[
\int _{0}^{\infty }t\, d\mu _{\tau }(t)=e^{a\tau },\quad \tau \ge 0.\]
There is a neighborhood $V$ of $0\in \Bbb C$ such that the map $\tau \mapsto \eta _{\mu _{\tau }}(z)$
is differentiable at $\tau =0$ for every $z\in \Omega \cup V$, and
the derivative \[
A(z)=\left.\frac{d\eta _{\mu _{\tau }}(z)}{d\tau }\right|_{\tau =0}\]
 is an analytic function of $z$. Moreover, we can write $A(z)=z(B(z)+a)$,
where $B$ is analytic in $\Omega \cup V$, $B(0)=0$, $B(\overline{z})=\overline{B(z)}$
and $\Im B(z)\ge 0$ whenever $\Im z>0$.

Conversely, for any $a\in \Bbb R$, and any analytic function $B$
defined in a set of the form $\Omega \cup V$, with $V$ a neighborhood
of $0$, satisfying the conditions above, there exists a unique $\circlearrowright _{0}$-semigroup
$\{\mu _{\tau }:\tau \ge 0\}$ of compactly supported probability
measures on $\Bbb R_{+}$such that\[
\left.\frac{d\eta _{\mu _{\tau }}(z)}{d\tau }\right|_{\tau =0}=z(B(z)+a),\quad z\in \Omega ,\]
and $\int _{0}^{\infty }t\, d\mu _{\tau }(t)=e^{a\tau }$ for $\tau \ge 0$.
Moreover, $\eta _{\mu _{t}}(z)=u_{\tau }(e^{a\tau }z)$, where $u_{\tau }(z)$
is the solution of the initial value problem\[
\frac{du_{\tau }(z)}{d\tau }=u_{\tau }(z)B(u_{\tau }(z)),\quad u_{0}(z)=z\in \Omega .\]
This solution exists for all $\tau \ge 0$.
\end{thm}
\begin{proof}
Start first with a semigroup $\{\mu _{\tau }:\tau \ge 0\}$, and define
functions $u_{\tau }:\Omega \to \Omega $ by setting $u_{\tau }(z)=\eta _{\mu _{\tau }}(e^{-a\tau }z)$
for $z\in \Omega $. Clearly then $u_{\tau }(z)$ depends continuously
on $z$, and the semigroup property can be translated into $u_{\tau +\tau '}(z)=u_{\tau }(u_{\tau '}(z))$.
As shown by Berkson and Porta {[}3{]} (see Theorem 1.1), these conditions
imply that $u_{\tau }(z)$ is a differentiable function of $\tau $,
and it satisfies the equation \[
\frac{du_{\tau }(z)}{d\tau }=C(u_{\tau }(z)),\]
where $C(z)=(du_{\tau }(z)/d\tau )|_{\tau =0}$. The initial condition
$u_{0}(0)=z$ comes from the identity $u_{0}=\eta _{\mu _{0}}=\eta _{\delta _{1}}$,
and this last function is easily seen to be the identity function
on $\Omega $. Clearly\[
\left.\frac{d\eta _{\mu _{\tau }}(z)}{d\tau }\right|_{\tau =0}=C(z)+az,\quad z\in \Omega .\]

Let us observe next that $\tau \mapsto \arg u_{\tau }(z)$ is an increasing
function for $\Im z>0$, and therefore\[
\Im \frac{C(z)}{z}=\left.\frac{d\log u_{\tau }(z)}{d\tau }\right|_{\tau =0}\ge 0,\]
so that indeed $C(z)=zB(z)$, where $B$ is an analytic function with
positive imaginary part in the upper half-plane $\Bbb C^{+}$. Moreover,
the fact that $u_{\tau }(0)=0$ yields $C(0)=0$, so that $B$ is
analytic in a neighborhood of zero. We also have $u'_{\tau }(0)=1$,
which shows that $C$ also has zero derivative at $z=0$, and therefore
$B(0)=0$ as well.

Conversely, assume that we are given an analytic function $B$ in
$\Omega \cup V$, with $B(0)=0$, and with positive imaginary part
in $\Bbb C^{+}$. We show first that the initial value problem\begin{equation}
\frac{du_{\tau }(z)}{d\tau }=u_{\tau }(z)B(u_{\tau }(z)),\quad u_{\tau }(0)=z\in \Omega \label{eq:(1)}\end{equation}
has a solution defined for all positive $\tau $. In order to do this
we apply another result of {[}3{]} (see Theorem 2.6, and the description
of the class $\mathcal{G}_{2}(\mathcal{H})$ for $b=0)$, which we
reformulate for the upper half-plane $\Bbb C^{+}$ and the left half-plane
$i\Bbb C^{+}$: Let $C:\Bbb C^{+}\to \Bbb C$ (resp., $C:i\Bbb C^{+}\to \Bbb C$)
be an analytic function such that $C(z)/z^{2}\in \Bbb C^{+}$ (resp.,
$-C(z)/z^{2}\in i\Bbb C^{+}$) for every $z$. Then for every $z$
(in the relevant domain), the initial value problem $du_{\tau }(z)/d\tau =C(u_{\tau }(z))$,
$u_{0}(z)=z$, has a solution defined for all positive $\tau $. The
function $C(z)=zB(z)$ satisfies the hypotheses of both of these results.
Indeed, the fact that $B$ has positive imaginary part in $\Bbb C^{+}$
allows us to write $B$ in Nevanlinna form\[
B(z)=\beta +\gamma z+\int _{-\infty }^{\infty }\frac{1+zt}{t-z}\, d\rho (t),\quad z\in \Bbb C^{+},\]
where $\beta $ is a real number, $\gamma \ge 0$, and $\rho $ is
a finite, positive Borel measure on $\Bbb R$. The fact that $B$
is real and analytic on $(-\infty ,\varepsilon ]$ for some $\varepsilon >0$
shows that the support of $\rho $ is contained in $[\varepsilon ,+\infty )$,
and the condition $B(0)=0$ yields the value\[
\beta =-\int _{0}^{\infty }\frac{1}{t}\, dt.\]
We conclude that\[
C(z)=z^{2}\left(\gamma +\int _{0}^{\infty }\frac{t^{2}+1}{t(t-z)}\, d\rho (t)\right),\quad z\in \Omega .\]
It is now easy to see that the integral above has positive imaginary
part if $z\in \Bbb C^{+}$, and positive real part for $z\in i\Bbb C^{+}$.
We conclude that the equation (3.1) has a solution defined for $\tau \ge 0$
for initial values $z$ in $\Bbb C^{+}\cup i\Bbb C^{+}$, and by symmetry
for all $z\in \Omega $. This equation will also have a solution defined
for small ,$\tau $ given an initial value $z>0$ sufficiently close
to zero. We deduce that, for small values of $\tau ,$ the function
$u_{\tau }$ is also analytic in a neighborhood of zero. The equation
$u_{\tau +\tau '}=u_{\tau }\circ u_{\tau '}$ shows that the same
is true for all values of $\tau $, and $u_{\tau }(0)=0$. The fact
that $B$ has positive imaginary part in $\Bbb C^{+}$ implies that
the function $\tau \mapsto \arg u_{\tau }(z)$ is an increasing function
of $\tau $, and therefore $\arg u_{\tau }(z)\ge \arg z$ for $z\in \Bbb C^{+}$.
(Note that $u_{\tau }(\Bbb C^{+})\subset \Bbb C^{+}$ by the theorem
of Berkson and Porta.) We conclude that there exist compactly supported
probability measures $\mu _{\tau }$ on $\Bbb R_{+}$ such that $\eta _{\mu _{\tau }}(z)=u_{\tau }(e^{at}z)$
for $z\in \Omega $ and $\tau \ge 0$. It is easy to verify now that
these measures form a multiplicative monotone convolution semigroup
satisfying the required conditions. The uniqueness of the semigroup
is a consequence of the uniqueness of solutions of ordinary differential
equations with a locally Lipschitz right-hand side.
\end{proof}
The results of Berkson and Porta {[}3{]} can also be formulated, via
conformal map, for the entire region $\Omega $. The corresponding
formulation however does not reflect the additional symmetries present
in our particular case.

The representation of the function $C$ found in the preceding proof
provides a bijection between $\circlearrowright _{0}$-convolution
semigroups and triples $(\gamma ,\rho ,a)$, where $a$ is a real
number, $\gamma \ge 0$, and $\rho $ is a finite, positive Borel
measure on some interval $[\varepsilon ,+\infty )$. The representation
of the function $A$ can be written more compactly if we use the measure
$\nu $ defined on the interval $[0,1/\varepsilon ]$ by the requirements
that $\nu (\{0\})=\gamma $ and $d\nu (t)=(t^{2}+1)d\rho (1/t)$ on
$(0,1/\varepsilon ]$. We have then\[
A(z)=az+z^{2}\int _{0}^{\infty }\frac{1}{1-zt}d\nu (t),\]
with $a\in \Bbb R$ and $\nu $ a positive, Borel, compactly supported
measure on $\Bbb R_{+}$. The constant $a$ is equal to zero if the
measures $\mu _{\tau }$ have first moment equal to one, in which
case the functions $\eta _{\mu _{\tau }}=u_{\tau }$ simply form a
semigroup relative to composition of functions on $\Omega $.

It is difficult to find explicit formulas for these semigroups. One
case when this is possible is $A(z)=\gamma z^{2}$ for some $\gamma >0$.
In this case the differential equation is easily solved, and it yields\[
\eta _{\mu _{\tau }}(z)=\frac{z}{1-\gamma \tau z},\quad \psi _{\mu _{\tau }}(z)=\frac{z}{1-(1+\gamma \tau )z},\quad z\in \Omega ,\]
so that\[
\mu _{\tau }=\frac{\gamma \tau }{1+\gamma \tau }\delta _{0}+\frac{1}{1+\gamma \tau }\delta _{1+\gamma \tau },\quad \tau \ge 0.\]

As in the case of additive monotone convolution {[}6{]}, the preceding
parametrization of semigroups also yields a parametrization of $\circlearrowright _{0}$-infinitely
divisible measures. Naturally, a compactly supported probability measure
$\mu $ on $\Bbb R_{+}$ is said to be $\circlearrowright _{0}$-\emph{infinitely
divisible} if, for every positive integer $n$, there exists a compactly
supported probability measure $\mu _{\frac{1}{n}}$ on $\Bbb R_{+}$
such that\[
\mu =\underbrace{\mu _{\frac{1}{n}}\circlearrowright _{0}\mu _{\frac{1}{n}}\circlearrowright _{0}\cdots \circlearrowright _{0}\mu _{\frac{1}{n}}}_{n\, \, \text {times}}.\]

\begin{thm}
Let $\mu \ne \delta _{0}$ be a $\circlearrowright _{0}$-infinitely
divisible, compactly supported, probability measure on $\Bbb R_{+}$.
There exists a unique $\circlearrowright $-semigroup $\{\mu _{\tau }:\tau \ge 0\}$
of compactly supported probability measures on $\Bbb R_{+}$ such
that $\mu _{1}=\mu $.
\end{thm}
\begin{proof}
Replacing the measure $\mu $ by the measure $d\mu (t/b)$, with $b=\int _{0}^{\infty }t\: d\mu (t)$
allows us to restrict ourselves to measures with first moment equal
to one. In this case it is clear that the measures $\mu _{\frac{1}{n}}$
satisfy the same property, and \[
\eta _{\mu }=\underbrace{{\eta _{\mu _{\frac{1}{n}}}\circ \eta _{\mu _{\frac{1}{n}}}\circ \cdots \circ \eta _{\mu _{\frac{1}{n}}}}}_{n\, \, \text {times}}.\]
The argument of Proposition 5.4 in {[}6{]} shows then that the measures
$\mu _{\frac{1}{n}}$are uniquely determined, so that we can further
define\[
\mu _{\frac{m}{n}}=\underbrace{\mu _{\frac{1}{n}}\circlearrowright _{0}\mu _{\frac{1}{n}}\circlearrowright _{0}\cdots \circlearrowright _{0}\mu _{\frac{1}{n}}}_{m\, \, \text {times}}\]
for arbitrary positive integers $m,n$. Clearly we have $\mu _{\tau +\tau '}=\mu _{\tau }\circlearrowright _{0}\mu _{\tau '}$
for rational $\tau ,\tau '>0$. It is then seen from Proposition 3.5
that the measures $\mu _{\frac{m}{n}}$ have uniformly bounded supports
if $m/n$ varies in a bounded set of rational numbers. We can now
verify that the map $\tau \mapsto \mu _{\tau }$ is continuous for
$\tau >0$ rational. Assume indeed that $\nu $ is the weak limit
of a sequence $\mu _{\tau _{k}}$, where $\tau _{k}>0$ are rational
numbers such that $\tau _{k}\to \frac{m}{n}$ as $k\to \infty $.
The continuity of multiplicative monotone convolution implies that
\[
\underbrace{{\nu \circlearrowright _{0}\nu \circlearrowright _{0}\cdots \circlearrowright _{0}\nu }}_{n\, \, \text {times}}=\mu _{m},\]
and the uniqueness of roots gives then $\nu =\mu _{\frac{m}{n}}$.
On the other hand, if $\tau _{k}\to 0$ and $\mu _{\tau _{k}}$ tends
to $\nu $, the measure $\nu \circlearrowright _{0}\mu $ is the weak
limit of $\mu _{1+\tau _{k}}$, so that $\nu \circlearrowright _{0}\mu =\mu $.
In other words, $\eta _{\nu }\circ \eta _{\mu }=\eta _{\mu }$, which
shows that $\eta _{\nu }$ must be the identity function, and hence
$\mu =\delta _{1}$. It is now easy to see that $\eta _{\tau }$ can
be defined for arbitrary $\tau >0$ by continuity. Indeed, consider
two sequences of positive rational numbers $\tau _{k}\to \tau ,\tau '_{k}\to \tau $
such that the sequences $\mu _{\tau _{k}},\mu _{\tau '_{k}}$ tend
weakly to measures $\nu ,\nu '$. By dropping to subsequences (and
possibly switching the two sequences), we may assume that $\tau _{k}>\tau '_{k}$
for all $k$. Since $\mu _{\tau _{k}}=\mu _{\tau '_{k}}\circlearrowright _{0}\mu _{\tau _{k}-\tau '_{k}}$,
and $\tau _{k}-\tau '_{k}\to 0$, we deduce that $\nu =\nu '$. The
uniqueness of the semigroup obtained this way follows immediately
from the uniqueness of $\mu _{\frac{1}{n}}$.
\end{proof}
There are analogous results for $\circlearrowright $-semigroups.

\begin{thm}
Consider a $\circlearrowright $-semigroup $\{\mu _{\tau }:\tau \ge 0\}$
of compactly supported measures on $\Bbb R_{+}$. The map $\tau \mapsto \eta _{\mu _{\tau }}(z)$
is differentiable for every $z\in \Omega $, and\[
\frac{d\eta _{\mu _{\tau }}(z)}{d\tau }=A(\eta _{\mu _{\tau }}(z)),\quad \tau \ge 0,z\in \Omega ,\]
where\[
A(z)=\left .\frac{d\eta _{\mu _{\tau }}(z)}{d\tau }\right |_{\tau =0},\quad z\in \Omega .\]
The function $A$ can be written as $A(z)=zB(z)$, where $B$ is analytic
in $\Omega $ and in a neighborhood of zero, and $\Im B(z)\ge 0$
for $z\in \Bbb C^{+}$.

Conversely, if $A$ is an analytic function in $\Omega $ with the
above properties, there exists a unique $\circlearrowright $-semigroup
$\{\mu _{\tau }:\tau \ge 0\}$ of compactly supported measures on
$\Bbb R_{+}$such that $A(z)=d\eta _{\mu _{\tau }}(z)/d\tau |_{\tau =0}$
for $z\in \Omega $.
\end{thm}
\begin{proof}
The differentiability of the map $\tau \mapsto \eta _{\mu _{\tau }}(z)$
follows from Theorem 1.1 of {[}3{]}, and the fact that $B$ has positive
imaginary part follows as before from the fact that the map $\tau \mapsto \arg \eta _{\mu _{\tau }}(z)$
is increasing when $z\in \Bbb C^{+}$. The uniqueness of the semigroup
$\mu _{\tau }$ is an immediate consequence of the uniqueness of solutions
to differential equations (with locally Lipschitz right-hand side).
The only thing that requires attention is the fact that, given a function
$A$ with the properties in the statement, the initial value problem\[
\frac{du}{d\tau }=A(u),\quad u(0)=z\in \Omega \]
has a solution defined for all $\tau \ge 0$. We will show that this
is in fact true whenever $B(z)=A(z)/z$ has positive imaginary part
in $\Bbb C^{+}$(without assuming that $B$ is analytic at zero).
To do this we write $B$ in Nevanlinna form\[
B(z)=\beta +\gamma z+\int _{0}^{\infty }\frac{1+zt}{t-z}\, d\rho (t),\quad z\in \Omega ,\]
with $\beta \in \Bbb R,$ $\gamma \in \Bbb R_{+}$, and $\rho $ a
positive Borel measure on $\Bbb R_{+}$. We will distinguish three
cases, according to the behavior of the function $B$ on the interval
$(-\infty ,0)$. Note that $B$ is increasing on this interval, so
that it could be negative on $(-\infty ,0)$, positive on $(-\infty ,0)$,
or vanish at some point in $(-\infty ,0)$. The first situation, $B(z)\le 0$
for all $z\in (-\infty ,0)$, amounts to $B(0-)\le 0$, which implies
that $\int _{0}^{\infty }\frac{1}{t}\, d\rho (t)$ is finite. After
rewriting the above formula as\[
B(z)=\beta +\int _{0}^{\infty }\frac{1}{t}\, d\rho (t)+\gamma z+z\int _{0}^{\infty }\frac{t^{2}+1}{t(t-z)}\, d\rho (t),\quad z\in \Omega ,\]
we deduce that\[
\beta +\int _{0}^{\infty }\frac{1}{t}\, d\rho (t)\le 0.\]
It is then easy to verify that $A(z)/z^{2}\in \Bbb C^{+}$ for $z\in \Bbb C^{+}$
and $-A(z)/z^{2}\in i\Bbb C^{+}$ for $z\in i\Bbb C^{+}$. We deduce
as in the proof of Theorem 3.6 that the solution to our initial value
problem extend to all $\tau \ge 0$. Assume next that $B(z)\ge 0$
for all $z\in (-\infty ,0)$. Since\[
\int _{0}^{\infty }\frac{1+zt}{t-z}\, d\rho (t)=o(z)\]
as $z\downarrow -\infty $, this is only possible when $\gamma =0$.
In this case\[
\lim _{z\downarrow -\infty }B(z)=\beta -\int _{0}^{\infty }t\, d\rho (t),\]
and we conclude that $\int _{0}^{\infty }t\, d\rho (t)<\infty $,
and $\beta \ge \int _{0}^{\infty }t\, d\rho (t)$. Setting $\alpha =\beta -\int _{0}^{\infty }t\, d\rho (t)$,
we have\[
A(z)=z\left (\alpha +\int _{0}^{\infty }\frac{t^{2}+1}{t-z}\, d\rho (t)\right ).\]
Using this formula, the inequality $\alpha \ge 0$, and the fact that
\[
\frac{z}{t-z}=-1+\frac{t}{t-z},\]
it is easy to see that $A(z)\in \Bbb C^{+}$ for $z\in \Bbb C^{+}$,
and $A(z)\in i\Bbb C^{+}$ for $z\in i\Bbb C^{+}$. The results of
Berkson and Porta show again that the solution $u$ of the initial
value problem extends to $\tau \ge 0$ for every $z\in \Bbb C^{+}\cup i\Bbb C^{+}$
and, by symmetry, for every $z\in \Omega $. (Note that in this case
the relevant Denjoy-Wolff point is infinity, which corresponds with
the family $\mathcal{G}_{1}(\mathcal{H})$ in the notation of {[}3{]}.)
Finally, assume that $B(-a)=0$ for some $a>0$. This yields the value\[
\beta =-\gamma a-\int _{0}^{\infty }\frac{1-at}{t+a}\, d\rho (t),\]
yielding the formula\[
A(z)=z(z+a)\left (\gamma +\int _{0}^{\infty }\frac{t^{2}+1}{(t+a)(t-z)}\, d\rho (t)\right ),\quad z\in \Omega .\]
 As in the preceding case, it will suffice to show that the initial
value problem for $u$ has a solution defined for all $\tau \ge 0$
if $z\in \Bbb C^{+}\cup i\Bbb C^{+}$. Using the results of {[}3{]}
(specifically, the classes $\mathcal{G}_{2}(\Bbb C^{+})$ and $\mathcal{G}_{3}(i\Bbb C^{+})$),
we see that $A$ must satisfy the following conditions:\[
\frac{A(z)}{(z+a)^{2}}\in \Bbb C^{+}\quad \text {for}\quad z\in \Bbb C^{+},\]
and\[
\frac{A(z)}{(z+a)(z-a)}\in -i\Bbb C^{+}\quad \text {for}\quad z\in i\Bbb C^{+}.\]
For the first of these conditions we write\[
\frac{A(z)}{(z+a)^{2}}=\frac{\gamma z}{z+a}+\int _{0}^{\infty }\frac{t^{2}+1}{t+a}\cdot \frac{z}{(t-z)(z+a)}\, d\rho (t),\]
which allows the calculation of the imaginary part\[
\Im \frac{A(z)}{(z+a)^{2}}=\Im z\left (\frac{\gamma a}{|z+a|^{2}}+\int _{0}^{\infty }\frac{t^{2}+1}{t+a}\cdot \frac{ta+|z|^{2}}{|t-z|^{2}|z+a|^{2}}\, d\rho (t)\right ).\]
This is clearly positive for $z\in \Bbb C^{+}.$ For the second condition
we have\[
\frac{A(z)}{(z+a)(z-a)}=\frac{\gamma z}{z-a}+\int _{0}^{\infty }\frac{t^{2}+1}{t+a}\cdot \frac{z}{(t-z)(z-a)}\, d\rho (t),\]
and\[
\Re \frac{A(z)}{(z+a)(z-a)}=\frac{\gamma (|z|^{2}-a\Re z)}{|z-a|^{2}}+\int _{0}^{\infty }\frac{t^{2}+1}{t+a}\cdot \frac{t|z|^{2}+a|z|^{2}-(ta+|z|^{2})\Re z}{(t-z)(z-a)}\, d\rho (t).\]
This is clearly positive when $\Re z<0$.

We have thus shown that the initial value problem has a solution defined
for all $\tau \ge 0$. Denote by $\eta _{\tau }(z)$ this solution.
This is an analytic function of $z$, and it extends analytically
to a neighborhood of zero if, in addition, $B$ is analytic at zero;
moreover, $\eta _{\tau }(0)=0$ in this case. It is shown now as in
the proof of Theorem 3.6 that $\eta _{t}=\eta _{\mu _{\tau }}$ for
some compactly supported measure $\mu _{\tau }$ on $\Bbb R^{+}$,
and these measures form a $\circlearrowright $-semigroup.
\end{proof}
As in the case of the operation $\circlearrowright $, $\delta _{0}$
is $\circlearrowright $-infinitely divisible. All other $\circlearrowright $-infinitely
divisible measures belong to a $\circlearrowright $-semigroup.

\begin{thm}
Let $\mu \ne \delta _{0}$ be a $\circlearrowright $-divisible measure,
compactly supported, probability measure on $\Bbb R^{+}$. There exists
a unique $\circlearrowright $-semigroup $\{\mu _{\tau }:\tau \ge 0\}$
of compactly supported probability measures on $\Bbb R_{+}$ such
that $\mu _{1}=\mu $.
\end{thm}
\begin{proof}
The argument is virtually identical with that of Theorem 3.7, except
that we need not start by normalizing the measures $\mu $. The details
are left to the interested reader.
\end{proof}

\section{Measures on the Unit Circle}

If $\mu $ is a probability measure on the unit circle $\Bbb T=\{\zeta \in \Bbb C:|\zeta |=1\}$,
the formal power series $\psi _{\mu },\eta _{\mu }$ converge in the
unit circle $\Bbb D=\{z\in \Bbb C:|z|<1\}$, and their sums are given
by\[
\psi _{\mu }(z)=\int _{\Bbb T}\frac{z\zeta }{1-z\zeta }\, d\mu (\zeta ),\quad \eta _{\mu }(z)=\frac{\psi _{\mu }(z)}{1+\psi _{\mu }(z)},\quad z\in \Bbb D.\]
An analytic function $\eta :\Bbb D\to \Bbb C$ is of the form $\eta _{\mu }$,
for some probability measure on $\Bbb T$, if and only if $|\eta (z)|\le |z|$
for all $z\in \Bbb T$ (cf., for instance, {[}1{]}). As in the case
of compactly supported measures on $\Bbb R_{+}$, the collection of
probability measures on $\Bbb T$ is identified with a subset of $\mathfrak{M}$.
The topology of $\mathfrak{M}$, restricted to this subset, is exactly
the topology of weak convergence of probability measures. One should
note that an element of $\mathfrak{M}$ may correspond to a measure
on $\Bbb T$, or to a measure on $\Bbb R_{+}$, and these two measures
can be quite different. The simplest occurence is the equality $\eta _{\delta _{0}}=\eta _{m}=0$,
where $\delta _{0}$ is a unit mass at the origin, while $m$ is normalized
arclength (or Haar) measure on $\Bbb T$.

\begin{prop}
If $\mu _{1},\mu _{2}$ are probability measures on $\Bbb T$ then
$\mu _{1}\circlearrowright \mu _{2}$ and $\mu _{1}\circlearrowright _{0}\mu _{2}$
are also probability measures on $\Bbb T$.
\end{prop}
\begin{proof}
If $|\alpha |\le 1$ is a complex number, we have\[
\left |\frac{1}{\alpha }\eta _{\mu _{2}}(\alpha z)\right |\le |z|,\quad z\in \Bbb D,\]
where the left-hand side must be interpreted as $|\eta '_{\mu _{2}}(0)z|$
when $\alpha =0$. We deduce that\[
\left |\eta _{\mu _{1}}\left (\frac{1}{\alpha }\eta _{\mu _{2}}(\alpha z)\right )\right |\le \left |\frac{1}{\alpha }\eta _{\mu _{2}}(\alpha z)\right |\le |z|,\quad z\in \Bbb D,\]
showing that the formal power series $\eta _{\mu _{1}\circlearrowright _{0}\mu _{2}}(z)$
corresponds indeed with a probability measure on $\Bbb T$. The measure
$\mu _{1}\circlearrowright \mu _{2}$ is treated similarly.
\end{proof}
The part of the preceding result concerning $\circlearrowright $
can be viewed as a consequence of the fact that the product of two
unitary operators is again a unitary operator. Indeed, given probability
measures $\mu _{1},\mu _{2}$ on $\Bbb T$, we can find unitary operators
$x_{1},x_{2}$ such that $x_{1}-1,x_{2}-1$ are monotonically independent,
and the distrinution of $x_{j}$ is $\mu _{j}$ for $j=1,2.$ It would
be nice to also understand the part concerning $\circlearrowright _{0}$
in the same manner, but it is not clear how to construct unitary operators
$x_{1},x_{2}$, with given distributions, such that $x_{1}-\varphi (x_{1}),x_{2}-\varphi (x_{2})$
are monotonically independent. Such operators are easily seen not
to exist in the standard realization used in Section 2.

Monotonic convolution semigroups of probability measures on $\Bbb T$
are defined as in the case of the half-line, and the following result
is the analogue of Theorem 3.6 in this context.

\begin{thm}
Consider a $\circlearrowright _{0}$-semigroup $\{\mu _{\tau }:\tau \ge 0\}$
of probability measures on $\Bbb T$. The map $\tau \mapsto \eta _{\mu _{\tau }}(z)$
is differentiable for every $z\in \Bbb D$, and the derivative \[
A(z)=\left.\frac{d\eta _{\mu _{\tau }}(z)}{d\tau }\right|_{\tau =0}\]
 is an analytic function of $z$. Moreover, we can write $A(z)=zB(z)$,
where $B$ is analytic in $\Bbb D$ and $\Re B(z)\le 0$ for $z\in \Bbb D$.

Conversely, for any analytic function $B$ defined in $\Bbb D$, with
$\Re B(z)\le 0$ for $z\in \Bbb D$, there exists a unique $\circlearrowright _{0}$-semigroup
$\{\mu _{\tau }:\tau \ge 0\}$ of probability measures on $\Bbb T$
such that\[
\left.\frac{d\eta _{\mu _{\tau }}(z)}{d\tau }\right|_{\tau =0}=zB(z),\quad z\in \Bbb D.\]
This semigroups satisfies $\int _{\Bbb T}\zeta \, d\mu _{\tau }(\zeta )=e^{B(0)\tau }$
for $\tau \ge 0$. Moreover, $\eta _{\mu _{t}}(z)=u_{\tau }(e^{B(0)\tau }z)$,
where $u_{\tau }:e^{B(0)\tau }\Bbb D\to \Bbb D$ is an analytic functions
satisfying the initial value problem\[
\frac{du_{t}(z)}{dt}=u_{t}(z)(B(u_{t}(z))-B(0)),\quad u_{0}(z)=z\in e^{B(0)\tau }\Bbb D.\]
This solution exists and belongs to $\Bbb D$ for all $t\in [0,\tau ]$.
\end{thm}
\begin{proof}
The numbers $\alpha (\tau )=\int _{\Bbb T}\zeta \, d\mu _{\tau }(\zeta )$
depend continuously on $\tau $, $\alpha (\tau +\tau ')=\alpha (\tau )\alpha (\tau ')$,
and $|\alpha (\tau )|\le 1$ for all $\tau $. It follows that $\alpha (\tau )=e^{a\tau }$
for some complex number $a$ with $\Re a\le 0$. Define now functions
$u_{\tau }:e^{a\tau }\Bbb D\to \Bbb D$ by $u_{\tau }(z)=\eta _{\mu _{\tau }}(e^{-a\tau }z)$
for $z\in e^{a\tau }\Bbb D$. These functions are analytic, and they
satisfy the equation\[
u_{\tau }(u_{\tau '}(z))=u_{\tau +\tau '}(z),\quad z\in e^{a(\tau +\tau ')}\Bbb D.\]
Moreover, tha map $t\mapsto u_{t}(z)$ is easily seen to be continuous
on the interval $[0,\tau ]$, provided that $z\in e^{a\tau }\Bbb D$.
The argument in Theorem 1.1 of {[}3{]} applies in this situation as
well, and it implies that the map $t\mapsto u_{t}(z)$ is in fact
differentiable, and the function\[
F(z)=\left .\frac{du_{\tau }(z)}{d\tau }\right |_{0},\quad z\in \Bbb D\]
is analytic. It follows that the map $\tau \mapsto \eta _{\mu _{\tau }}(z)$
is differentiable as well, and the function $A$ in the statement
is analytic. In fact, we have $A(z)=F(z)-az$ since $\eta _{\mu _{0}}(z)=z$.
In order to show that $A$ has the required form, let us also consider
the function $v_{\tau }(z)=e^{a\tau }u_{\tau }(z)=e^{a\tau }\eta _{\mu _{\tau }}(e^{-a\tau }z)$
defined in $e^{a\tau }\Bbb D$, for which\[
\left .\frac{dv_{\tau }(z)}{d\tau }\right |_{0}=az+\left .\frac{du_{\tau }(z)}{d\tau }\right |_{0}=A(z),\quad z\in \Bbb D.\]
For this function we have $|v_{\tau }(z)|\le |z|=|v_{0}(z)|,$ so
that indeed\[
\Re \frac{A(z)}{z}=\left .\frac{d\Re \log v_{\tau }(z)}{d\tau }\right |_{\tau =0}=\left .\frac{d\log |v_{\tau }(z)|}{\tau }\right |_{\tau =0}\le 0,\quad z\in \Bbb D\setminus \{0\}.\]
Let us then write $A(z)=zB(z)$, and verify that $a=-B(0)$. Indeed,
all the functions $(u_{\tau }(z)-z)/\tau $ have a double zero at
the origin, and therefore so does their limit $F(z)$; therefore $B(z)+a$
must be zero for $z=0$.

Conversely, assume that \textbf{$B$} is an analytic function with
negative real part in $\Bbb D$. It will suffice to show that the
initial value problem\[
\frac{du_{t}(z)}{dt}=u_{t}(z)(B(u_{t}(z))-B(0)),\quad u_{t}(0)=z\in e^{B(0)\tau }\Bbb D\]
has a solution defined on the entire interval $[0,\tau ]$, and that
\[
|u_{\tau }(z)|\le e^{-\Re B(0)\tau }|z|,\quad z\in e^{B(0)\tau }\Bbb D.\]
Indeed, once this is done, we can define the functions $\eta _{\tau }:\Bbb D\to \Bbb D$
by $\eta _{\tau }(z)=u_{\tau }(e^{B(0)\tau }z)$, and these functions
will be of the form $\eta _{\tau }=\eta _{\mu _{\tau }}$ for some
probability measures $\mu _{\tau }$ which are easily seen to form
a $\circlearrowright _{0}$-semigroup. The existence of the solutions
$u_{t}$ on the stated interval is easy to deduce from the general
theory of ordinary differential equations. We sketch a somewhat more
direct argument based on an appropriate approximation scheme. Namely,
define functions $w_{\varepsilon }:\Bbb D\to \Bbb C$ by\[
w_{\varepsilon }(z)=ze^{\varepsilon (B(z)-B(0))}\quad z\in \Bbb D,\varepsilon >0.\]
These functions satisfy $|w_{\varepsilon }(z)|\le e^{-\varepsilon B(0)}|z|$.
We then define $u_{\tau }^{(n)}:e^{B(0)\tau }\Bbb D\to \Bbb D$ by\[
u_{\tau }^{(n)}=\underbrace{w_{\tau /n}\circ w_{\tau /n}\circ \cdots \circ w_{\tau /n}}_{n\, \, \text {times}};\]
it is easy to see that $u_{\tau }^{(n)}$ is indeed defined in $e^{B(0)\tau }\Bbb D$.
There exists a positive number $\delta $ such that $u_{\tau }^{(n)}|\delta \Bbb D$
converge uniformly as $n\to \infty $ to the solution $u_{\tau }$
of our initial value problem, provided that $\tau \le \delta $. Now,
the functions $u_{\tau }^{(n)}$ are analytic and uniformly bounded
on $e^{B(0)\delta }\Bbb D$ for $\tau \le \delta $, and therefore
$\lim _{n\to \infty }u_{\tau }^{(n)}$ will exist (by the Vitali-Montel
theorem) on the entire disk $e^{B(0)\delta }\Bbb D$ for all such
$\tau $. In an analogous fashion, we deduce that $u_{\tau }(z)=\lim _{n\to \infty }u_{\tau }^{(n)}(z)$
exists for all $z\in e^{B(0)\tau }\Bbb D$ if $\tau \le \delta $.
Observe now the equality \[
u_{\tau }^{(n)}\circ u_{\tau '}^{(n')}=u_{\tau +\tau '}^{(n+n')}\quad \text {when}\quad \frac{\tau }{n}=\frac{\tau '}{n'},\]
which shows now that the convergence of $u_{\tau }^{(n)}$ can be
extended from the interval $[0,\delta ]$ to arbitrary $\tau >0$,
yielding a function $u_{\tau }$ defined in the common domain of $u_{\tau }^{(n)}$.
Clearly these functions will solve the initial value problem in the
required range. 
\end{proof}
The preceding result yields a parametrization of all $\circlearrowright _{0}$-semigroups
on the unit circle. In fact, every analytic function $B$ with negative
real part on $\Bbb D$ can be written using the Herglotz formula\[
B(z)=i\beta -\int _{\Bbb T}\frac{\zeta +z}{\zeta -z}\, d\rho (\zeta ),\quad z\in \Bbb D,\]
where $\beta $ is a real number, and $\rho $ is a finite positive
Borel measure on $\Bbb T$. The constant $a=B(0)$ is then given by\[
a=i\beta -\rho (\Bbb T),\]
and the differential equation for $u_{\tau }$ is\[
\frac{du_{t}(z)}{dt}=2u_{t}(z)^{2}\int _{\Bbb T}\frac{d\rho (\zeta )}{u_{t}(z)-\zeta },\quad u_{0}(z)=z\in e^{a\tau }\Bbb D.\]
As in the case of the half-line, the solutions of this equation can
seldom be calculated explicitly. The case $\rho =0$ corresponds with
semigroups where each $\mu _{\tau }$ is a point mass. In all cases
when $\rho \ne 0$, it is easy to see that the measures $\mu _{\tau }$
converge weakly to Haar measure $m$ as $\tau \to \infty $. One semigroup
which can be calculated explicitly corresponds with $B(z)=z^{n}-1$,
where $n\ge 1$ is an integer. We just mention the following formula:\[
u_{\tau }(z)=\frac{z}{(1-(n+1)z^{n}\tau )^{1/n}},\quad z\in e^{-\tau }\Bbb D,\]
where the root is chosen to be equal to one at the origin.

Infinite divisibility can also be characterized in terms of semigroups
in the case of the circle. As for the half-line (where $\delta _{0}$
$\circlearrowright _{0}$-infinitely divisible, but not part of a
semigroup), there is an exception, namely Haar measure $m$ which
satisfies $m\circlearrowright _{0}m=m\circlearrowright m=m$. More
generally, we have the following result.

\begin{lem}
If $\mu _{1},\mu _{2}$ are probability measures on $\Bbb T$, and
$\int _{\Bbb T}\zeta \, d\mu _{1}(\zeta )=\int _{\Bbb T}\zeta \, d\mu _{2}(\zeta )=0$,
then $\mu _{1}\circlearrowright _{0}\mu _{2}=m$.
\end{lem}
\begin{proof}
We have $\eta _{\mu _{1}\circlearrowright \mu _{2}}(z)=\eta _{\mu _{1}}(\eta '_{\mu _{2}}(0)z)=\eta _{\mu _{1}}(0)=0$
since $\eta '_{\mu _{1}}(0)=\eta '_{\mu _{2}}(0)=0$. Alternatively,
one observes that two monotonically independent variables $x_{1},x_{2}$
such that $\varphi (x_{1})=\varphi (x_{2})=0$ must satisfy $\varphi ((x_{1}x_{2})^{n})=0$
for all $n\ge 1$.
\end{proof}
We conclude that a $\circlearrowright _{0}$-infinitely divisible
probability measure $\mu $ on $\Bbb T$ with first moment zero must
in fact coincide with $m$. Indeed, $\mu =\mu _{\frac{1}{2}}\circlearrowright _{0}\mu _{\frac{1}{2}}$,
and the measure $\mu _{\frac{1}{2}}$ must also have first moment
equal to zero.

\begin{thm}
Let $\mu \ne m$ be a $\circlearrowright _{0}$-infinitely divisible
probability measure on $\Bbb T$. There exists a $\circlearrowright _{0}$-semigroup
$\{\mu _{\tau }:\tau \ge 0\}$ of probability measures on $\Bbb T$
such that $\mu _{1}=\mu $.
\end{thm}
\begin{proof}
As noted before the statement, we can write $\int _{\Bbb T}\zeta \, d\mu (\zeta )=\rho e^{i\theta }$
with $\theta \in \Bbb R$ and $\rho >0$. Choose for each interger
$n\ge 1$ a measure $\nu _{n}$ such that $\mu =\nu _{n}^{\circlearrowright _{0}2^{n}}$;
these measures are no longer uniquely determined, but (possibly after
an appropriate rotation) can be assumed to satisfy $\int _{\Bbb T}\zeta \, d\nu _{n}(\zeta )=\rho ^{1/2^{n}}e^{i\theta /2^{n}}$.
There exists a sequence $n_{1}<n_{2}<\cdots $ with the property that
the each sequence $\{\nu _{n_{j}}^{\circlearrowright _{0}2^{n_{j}-n}}:j\ge n\}$
has a weak limit; call this limit $\mu _{\frac{1}{2^{n}}}$. These
measures will then satisfy\[
\int _{\Bbb T}\zeta \, d\mu _{\frac{1}{2^{n}}}(\zeta )=\rho ^{1/2^{n}}e^{i\theta /2^{n}},\quad \mu _{\frac{1}{2^{n}}}^{\circlearrowright _{0}2^{n}}=\mu ,\quad \text {and}\quad \mu _{\frac{1}{2^{n}}}^{\circlearrowright _{0}2^{m}}=\mu _{\frac{1}{2^{n-m}}}\quad \text {for}\; m<n.\]
Note that the measures $\mu _{\frac{1}{2^{n}}}$ converge weakly to
$\delta _{1}$ as $n\to \infty $; indeed, their first moments converge
to $1$, and $\delta _{0}$ is the only probability measure on $\Bbb T$
with first moment equal to one. We can now define\[
\mu _{\frac{m}{2^{n}}}=\mu _{\frac{1}{2^{n}}}^{\circlearrowright _{0}m}\]
for $m,n$ positive integers, and this is a good definition, i.e.,
it depends only on the fraction $m/2^{n}$ and not on the value of
$n$. With this definition, it is still true that $\mu _{\tau }$
tends weakly to $\delta _{0}$ if $\tau \to 0$. Let now $\tau $
be an arbitrary positive number, and choose numbers $\tau _{k},\tau '_{k}$
of the form $m/2^{n}$ such that $\lim _{k\to \infty }\tau _{k}=\lim _{k\to \infty }\tau '_{k}=\tau $,
and the sequences $\{\mu _{\tau _{k}},k\ge 1\},\{\mu _{\tau '_{k}},k\ge 1\}$
have weak limits $\nu ,\nu '$. Dropping to subsequences we can assume
that $\tau _{k}<\tau '_{k}$ for all $k$. The equality $\mu _{\tau '_{k}}=\mu _{\tau '_{k}-\tau _{k}}\circlearrowright _{0}\mu _{\tau _{k}}$
yields then $\nu '=\delta _{0}\circlearrowright _{0}\nu =\nu $. This
unique limit can then be denoted $\mu _{\tau }$. It is easy to verify
that the measures $\mu _{\tau }$ form a multiplicative monotonic
convolution semigroup, and $\mu _{1}=\mu $.
\end{proof}
The semigroup provided by the preceding theorem is never unique. Thus,
if the semigroup is generated (in the sense of Theorem 4.2) by the
function $zB(z)$, then the function $z(B(z)+2\pi i)$ will generate
a new semigroup with $\mu _{1}=\mu $. Of course, the only difference
between these semigroups is a rotation of angle $2\pi \tau $ of the
measure $\mu _{\tau }$. It is fairly easy to see that this is the
only possible kind of nonuniqueness. More precisely, we have the following
result.

\begin{prop}
If $\mu ,\mu _{1},\mu _{2}\in \mathfrak{M}$ are such that $\mu _{1}\circlearrowright _{0}\mu _{1}=\mu _{2}\circlearrowright _{0}\mu _{2}=\mu $
and $\mu _{1}(X)=\mu _{2}(X)\ne 0$, then $\mu _{1}=\mu _{2}$. The
same result is true for the operation $\circlearrowright $.
\end{prop}
\begin{proof}
If $\mu _{1}(X)=\mu _{2}(X)=1$, then we have $\eta _{\mu _{1}}\circ \eta _{\mu _{1}}=\eta _{\mu _{2}}\circ \eta _{\mu _{2}}=\eta _{\mu }$.
In this case the result follows from the argument of Proposition 5.4
in {[}6{]}. The general case reduces to this particular one by considering
the new distributions $\nu _{j}(p(X))=\mu _{j}(p(X/\alpha ))$, $p\in \Bbb C[X]$,
where $\alpha =\mu _{1}(X)=\mu _{2}(X)$.
\end{proof}
This result shows that in fact the measures $\nu _{n}$ in the proof
of Theorem 4.4 are uniquely determined, and therefore there is precisely
one semigroup for every choice of the argument of $\int _{\Bbb T}\zeta \, d\mu (\zeta )$.

The analogue of Theorem 4.2 for $\circlearrowright $-semigroups is
obtained directly from the results of Berkson and Porta {[}3{]}. Indeed,
the corresponding functions $\eta _{\mu _{\tau }}$ simply form a
composition semigroup of analytic maps of the disk, fixing the origin.
We record the result below.

\begin{thm}
Consider a $\circlearrowright $-semigroup $\{\mu _{\tau }:\tau \ge 0\}$
of probability measures on $\Bbb T$. The map $\tau \mapsto \eta _{\mu _{\tau }}(z)$
is differentiable for every $z\in \Bbb D$, and the derivative \[
A(z)=\left.\frac{d\eta _{\mu _{\tau }}(z)}{d\tau }\right|_{\tau =0}\]
 is an analytic function of $z$. Moreover, we can write $A(z)=zB(z)$,
where $B$ is analytic in $\Bbb D$ and $\Re B(z)\le 0$ for $z\in \Bbb D$.

Conversely, for any analytic function $B$ defined in $\Bbb D$, with
$\Re B(z)\le 0$ for $z\in \Bbb D$, there exists a unique $\circlearrowright $-semigroup
$\{\mu _{\tau }:\tau \ge 0\}$ of probability measures on $\Bbb T$
such that\[
\left.\frac{d\eta _{\mu _{\tau }}(z)}{d\tau }\right|_{\tau =0}=zB(z),\quad z\in \Bbb D.\]
The functions $\eta _{\mu _{\tau }}$ satisfy the initial value problem\[
\frac{d\eta _{\mu _{\tau }}(z)}{d\tau }=\eta _{\mu _{\tau }}(z)B(\eta _{\mu _{\tau }}(z)),\quad \eta _{\mu _{\tau }}(0)=z\in \Bbb D.\]

\end{thm}
Infinite divisibility is also characterized in terms of semigroups,
and the remarks about uniqueness made about $\circlearrowright _{0}$-divisible
measures apply here as well. The proofs given above are easily converted
to this setting.

\begin{thm}
Let $\mu \ne m$ be a $\circlearrowright $-infinitely divisible probability
measure on $\Bbb T$. There exists a $\circlearrowright $-semigroup
$\{\mu _{\tau }:\tau \ge 0\}$ of probability measures on $\Bbb T$
such that $\mu _{1}=\mu $.
\end{thm}

\lyxaddress{Mathematics Department, Indiana University, Bloomington, IN 47405,
USA}
\end{document}